\documentclass[graybox]{svmult}

\usepackage{mathptmx}       
\usepackage{helvet}         
\usepackage{courier}        
\usepackage{type1cm}        
\usepackage{makeidx}         
\usepackage{graphicx}        
\usepackage{multicol}        
\usepackage[bottom]{footmisc}

\usepackage{amsmath,amsfonts,amssymb}
\usepackage{booktabs}
\usepackage{subfig}
\usepackage{algorithm,algpseudocode}
\newcommand{\ol}{\overline}
\newcommand{\ul}{\underline}
\newcommand{\dist}[1]{\underline{#1}_{\,dist}}
\newcommand{\glob}[1]{\underline{#1}_{\,glob}}
\newcommand{\acc}[1]{\underline{#1}_{\,acc}}
\newcommand{\Acc}{\mbox{$\Sigma$}}
\newcommand{\Accglob}{\Acc_{glob}\;}
\newcommand{\Split}{\mbox{\raisebox{.15em}{$\chi$}}_{glob}\;}

\makeindex             

\begin{document}

\title*{A parallel multigrid solver for multi-patch Isogeometric Analysis}
\author{Christoph Hofer and Stefan Takacs}
\institute{Christoph Hofer \at Doctoral Program Computational Mathematics, University Linz,
Altenberger Str. 69, 4040 Linz, \email{christoph.hofer@dkcm.jku.at}
\and Stefan Takacs \at RICAM, Austrian Academy of Sciences, Altenberger Str. 69, 4040 Linz,
\email{stefan.takacs@ricam.oeaw.ac.at}}
%
%
\maketitle

\abstract{
Isogeometric Analysis (IgA) is a framework for setting up spline-based discretizations
of partial differential equations, which has been introduced around a decade ago
and has gained much attention since then.
If large spline degrees are considered, one obtains the approximation power of a high-order
method, but the number of degrees of freedom behaves like for a low-order method.
One important ingredient to use a discretization with large spline degree, is a robust
and preferably parallelizable solver.
While numerical evidence shows that multigrid solvers with standard
smoothers (like Gauss Seidel) does not perform well if the spline degree
is increased, the multigrid solvers proposed by the authors and their co-workers
proved to behave optimal both in the grid size and the spline degree. In
the present paper, the authors want to show that those solvers are
parallelizable and that they scale well in a parallel environment.
}

\section{Introduction}
Isogeometric Analysis (IgA) was originally introduced in the seminal paper~\cite{Hughes:2005},
aiming to unite the worlds of computer aided design (CAD) and finite element (FEM) simulation.
From a technical point of view, it
is a framework for setting up spline-based discretizations
of partial differential equations. 
The key idea is that the spline space is typically first defined on the unit square or the unit
cube and then mapped to the computational domain using
one global geometry function. More complicated domains cannot be
represented by just one such geometry function. Instead, the 
computational domain is decomposed into patches, where each of them is
represented by its own geometry function. This is  called the
\emph{multi-patch case}, in contrast to the \emph{single-patch case}.

As a next step, the linear system resulting from the discretization
of the PDE has to be solved. This might be challenging as the
condition number of the linear system grows exponentially with the
spline degree, where high spline degrees might be desired because
of their superior approximation power.

While in early IgA literature, the dependence of methods on the
spline degree has not been considered, in the last few years
robustness in the spline degree has gained increasing interest.
Several (almost) robust approaches or approaches with a mild
dependence on the spline degree have been proposed, on the one side for the single-patch
case, cf.~\cite{CollierEtAl:2012,Donatelli:2014a,HTZ:2016,Sangalli:2016,HT:2016} and references therein, and
on the other side as approaches aiming to combine patch-local solvers
to a global solver,
cf.~\cite{KleissEtAl:2012,DaVeigaEtAl:2012,DaVeigaEtAl:2014,HLT:2017} and
references therein. 

In~\cite{Takacs:2017}, we have considered a slightly different approach: We
do not aim to combine patch-local solvers to a global solver, but to combine patch-local
smoothers to a global smoother which is used within a global multigrid
solver. In the present paper, we give some additional remarks on an efficient implementation
of the multigrid method, comment on its parallelization and give numerical results.

This paper is organized as follows. First, the model problem and
the discretization are discussed in Sec.~\ref{sec:prelim}. Then, in Sec.~\ref{sec:mg},
we recall the formulation of the multigrid solver. Its parallelization is discussed
in the following Sec.~\ref{sec:parallel}.  In Sec.~\ref{sec:num}, we give the results of
numerical experiments and draw conclusions.

\section{Model problem and isogeometric discretization}\label{sec:prelim}
 
Let $\Omega\subset \mathbb{R}^d$ with $d\in\{2,3\}$ be a bounded computational domain with Lipschitz boundary.
We consider a standard \emph{Poisson model problem} 
\[
		- \Delta u = f \quad\mbox{in}\quad \Omega\ , \quad
		u = 0 \quad\mbox{on}\quad \Gamma_D
		\qquad \mbox{and} \qquad
		\frac{\partial u}{\partial n} = 0 \quad\mbox{on}\quad \Gamma_N \ ,
\]
where $\Gamma_D$ is a subset of $\partial \Omega$ with positive measure and
$\Gamma_N:=\partial\Omega \backslash \Gamma_D$.
The model problem reads in variational form as follows. Given $f\in L_2(\Omega)$,
find $u\in H^1_{0,D}(\Omega)$ such that
\begin{equation} \label{eq:model}
			(\nabla u,\nabla v)_{L_2(\Omega)} = (f,v)_{L_2(\Omega)}
\qquad \mbox{for all} \quad v \in H^1_{0,D}(\Omega)\ .
\end{equation}
Here and in what follows, $L_2(\Omega)$ and $H^1(\Omega)$ are the standard Lebesgue and Sobolev spaces
with standard norms and
$H^1_{0,D}(\Omega):=\{u \in H^1(\Omega)\;:\; u|_{\Gamma_D} = 0\}$.

We preform a standard isogeometric multi-patch discretization as it has been specified
in~\cite{Takacs:2017}. In the present paper, we try to keep
the explanation short and give only an overview.
We assume that the computational domain $\Omega$ is composed of $K$ patches $\Omega_k$ such that
\begin{equation} \label{eq:matching}
	\ol{\Omega} = \bigcup_{k=1}^K \ol{\Omega_k}
	\quad \mbox{and} \quad
	\Omega_k\cap\Omega_l=\emptyset \mbox{ for any }k\not=l\ ,
\end{equation}
where each patch $\Omega_k$ is a bounded and open domain. We assume that the patches are
fully matching, i.e., the intersections
$\ol{\Omega_k}\cap\ol{\Omega_l}$ are either empty, common vertices, common edges or common faces.
Any of the patches is parametrized by a bijective geometry function
\[
		\textbf{G}_k :\widehat{\Omega}:=(0,1)^d \rightarrow \Omega_k
				:= \textbf{G}_k (\widehat{\Omega})\subset \mathbb{R}^d\ .
\]

\begin{figure}[h]
\begin{center}
	\includegraphics[width=.17\textwidth]{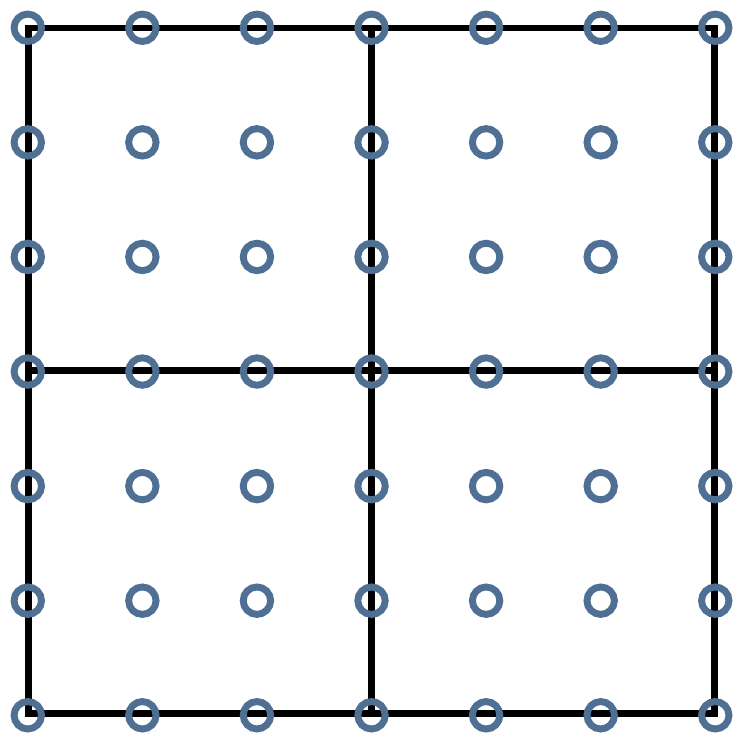}
\end{center}
\caption{Fully matching discretization}\label{fig:decomp1}
\end{figure}

Before we define set of trial functions $V_\ell \subset H^1_{0,D}(\Omega)$, we introduce discretizations
living on the parameter
domain~$\widehat{\Omega}$. Let
\[
	S_{p,h}(0,1) := \left\{ u \in C^{p-1}(0,1) \;:\; u|_{[hi,h(i+1)]} \mbox{ is a polynomial of degree } p \;,\; \forall_{i=1,\ldots,n} \right\}
\]
be the space of univariate splines of maximum smoothness and the space $S_{p,h}(\widehat{\Omega}):= 
S_{p,h}(0,1) \otimes \cdots \otimes S_{p,h}(0,1)$ be the corresponding tensor-product spline space.
The grid size $h$ and the spline degree $p$ might be different for any patch and for any spacial direction;
for simplicity, we do not express that in the notation.
Based on the discretization living on the parameter domain $\widehat{\Omega}$, we define the function
space $V_\ell$ of isogeometric functions living on the physical domain $\Omega$ as follows:
\begin{equation}\label{eq:vh}
	V_\ell := \{ u\in C^0(\Omega) \;:\; u\circ \textbf{G}_k \in S_{p,h_\ell}(\widehat{\Omega}) \} \ .
\end{equation}
We assume to have a \emph{fully matching discretization}, which means that the discretizations agree
on the interfaces. A more formal definition of the basis and its discretization is given in~\cite[Sec.~2]{Takacs:2017}. 
In Fig.~\ref{fig:decomp1}, a fully matching discretization is depicted, where each node represents one
basis function and therefore one degree of freedom (dof). Note that any of the basis function whose associated node lies
on one patch, vanishes outside of that patch. Any of the basis functions whose associated node lies within
one edge, vanishes outside the union of the edge and the adjacent patches. Finally, any of the basis
functions whose associated node coincides with one vertex, vanishes outside the union of that vertex and 
the adjacent edges and patches. The behavior in three dimensions is completely analogous.

The Galerkin principle yields the following discretized variational problem. Find $u\in V_\ell$ such that
\begin{equation} \label{eq:model:discr}
			a( u, v)= (f,v)_{L_2(\Omega)}
\qquad \mbox{for all $v \in V_\ell$}\ ,
\end{equation}
where 
\begin{equation} \label{eq:model:bil}
	a( u, v) := (\nabla u,\nabla v)_{L_2(\Omega)} = \sum_{k=1}^K  \underbrace{  ( |\det J_{\textbf{G}_k} | J_{\textbf{G}_k}^{-\top} J_{\textbf{G}_k}^{-1} \nabla \widehat{u}_k, \nabla \widehat{v}_k )_{L^2(\widehat{\Omega})} }_{\displaystyle a_k(u, v):= }
\end{equation}
for $\widehat{u}_k := u \circ \textbf{G}_k \in S_{p,h_\ell}(\widehat{\Omega})$ and
$\widehat{v}_k := v \circ \textbf{G}_k \in S_{p,h_\ell}(\widehat{\Omega})$
and where $J_{\textbf{G}_k}$ is the Jacobian of the geometry map.
Using the chosen basis, we obtain a matrix-vector formulation of the discretized problem,
which reads as follows. Find $\ul{u} \in \mathbb{R}^N$ such that
\begin{equation} \label{eq:linear:system}
			A_\ell \, \ul{u} = \ul{f}\ .
\end{equation}
Allowing constants that depend on the geometry function, we obtain that 
the matrix $A_\ell$ is spectrally equivalent to the matrix $\widehat{A}_\ell$, which discretizes the bilinear form
\[
	\widehat{a}\, (u, v):= \sum_{k=1}^K    (  \nabla \widehat{u}_k, \nabla \widehat{v}_k )_{L^2(\widehat{\Omega})} \ ,
\]
where, again, $\widehat{u}_k := u \circ \textbf{G}_k $ and $\widehat{v}_k := v \circ \textbf{G}_k $.

\section{The multigrid solver and its extension to three dimensions}\label{sec:mg}

We employ the multigrid solver based on a hierarchy of grids for grid levels $\ell=0,\ldots,L$, obtained by uniform refinement.
Throughout the grid hierarchy, the spline degree $p$ and the corresponding smoothness is kept
unchanged. This yields nested spaces: $V_0 \subset V_1 \subset \cdots \subset V_L \subset H^1_{0,D}(\Omega)$,
which allows to use the canonical embedding $V_{\ell-1}\rightarrow V_\ell$ for the multigrid method;
its matrix representation is denoted by $P_\ell$. Following the usual pattern, we use its transpose
$P_\ell^\top$ as restriction.

One \emph{multigrid cycle} on some grid level $\ell$ consists of the following steps.
\begin{itemize}
	\item First, $\nu$ \emph{pre-smoothing steps} are applied, where each reads as follows: 
			\begin{equation}\label{eq:sm}
				\ul{u} \gets \ul{u} + \tau L_\ell^{-1} (\ul{f} - A_\ell \ul{u}) \ .
			\end{equation}
		The choice of the smoothing operator $L_\ell^{-1}$ and the damping parameter $\tau$ are discussed below.
	\item Then, the \emph{coarse grid correction} is performed:
		\[
			\ul{u} \gets \ul{u} + \tau P_\ell A_{\ell-1}^{-1} P_\ell^\top (\ul{f} - A_\ell \ul{u}) \ ,
		\]
		where for $\ell>1$, the application $A_{\ell-1}^{-1}$ is replaced by $\mu=1$ (V-cycle) or 
		$\mu=2$ (W-cycle) recursive applications of the multigrid method on the coarser grid level.
	\item Finally, again $\nu$ \emph{post-smoothing steps} \eqref{eq:sm} are applied.
\end{itemize}

\begin{figure}[h]
\begin{center}
	\includegraphics[width=.17\textwidth]{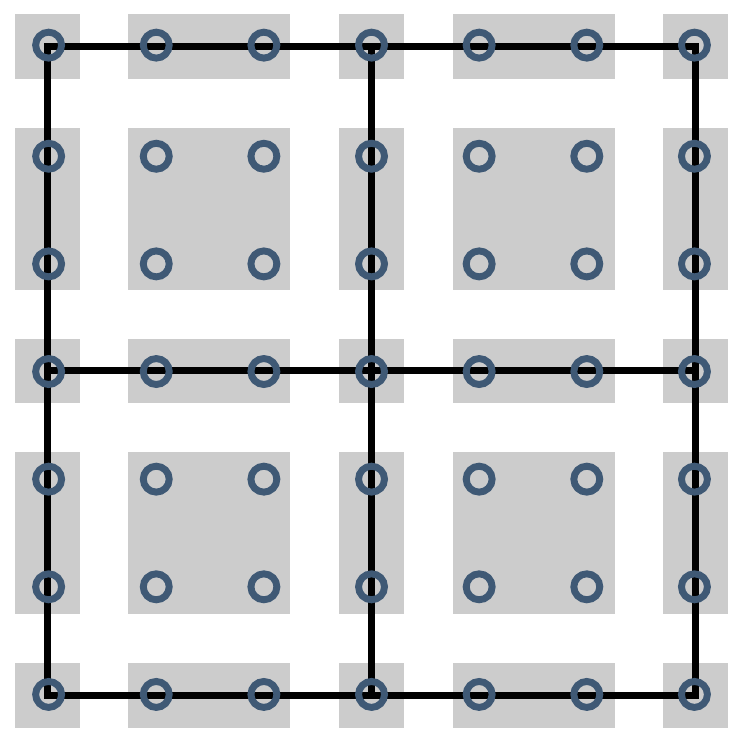}
\end{center}
\caption{Decomposition into pieces serving as subspaces for the additive Schwarz method}\label{fig:decomp2}
\end{figure}

As smoother, an additive Schwarz type combination 
\begin{equation}\label{eq:additive}
	L_\ell^{-1} := \sum_{T} P_{\ell,T} L_{\ell,T}^{-1} P_{\ell,T}^{\top}
\end{equation}
of local smoothing operators $L_{\ell,T}^{-1}$ is proposed, where the dofs are collected based on separating
the domain into \emph{pieces}: patches, vertices, edges and, in three dimensions, faces.
Here, each dof is assigned to exactly one of these pieces, cf. Fig.~\ref{fig:decomp2}.
Certainly, based on such a one-by-one splitting, the matrix $P_{\ell,T}$ is nothing but a indicator
matrix representing the canonical embedding.
\mbox{}\\
The local smoothing operators are chosen as follows.
\begin{itemize}
	\item For the patch-interiors, the subspace corrected mass smoother as proposed in~\cite{HT:2016} is
				chosen as smoothing operator $L_{\ell,T}^{-1}$.
	\item For the edges and vertices, in~\cite{Takacs:2017} direct solvers have been
				proposed as smoothers, i.e., $L_{\ell,T}$ is 
				the restriction of the matrix $A_\ell$ to the edge or vertex. To avoid unnecessary
				communication, we choose an approximation which can be computed directly.
				Using \cite[Lemma~4.1]{Takacs:2017} and \cite[eq.~(4.16)]{Takacs:2017}, 
				we obtain that the restriction of $A_\ell$ to an edge is spectrally equivalent to
				\[
					L_{\ell,T}:=
 					\left( \frac{h_\ell}p \right)^{d-1} \mathrm{K}_\ell
					+ 
					\left( \frac {h_\ell}p \right)^{d-3} \mathrm{M}_\ell \ ,
				\]
				where $\mathrm{K}_\ell$ and $\mathrm{M}_\ell$ are the corresponding univariate
				stiffness and mass matrices.
				Analogously, its restriction to a vertex is a constant in the order of
				\[
					L_{\ell,T}:=
						\left( \frac{h_\ell}p \right)^{d-2} \ .
				\]
	\item Three dimensional problems have not been considered in~\cite{Takacs:2017}, so we have to discuss how to
				choose the local smoothers for faces. If, as for the edges and vertices, again a direct solver was applied,
				the overall computational costs would not be optimal anymore. So, again, observe that the
				the restriction of $A_\ell$ to a face is spectrally equivalent to
				\[
					L_{\ell,T}^*:=
					 \left( \frac{h_\ell}p \right)^{d-2} \mathcal{K}_{\ell} + \left(  \frac {h_\ell} p \right)^{d-4} \mathcal{M}_{\ell}\ ,
				\]
				where and $\mathcal{K}_{\ell}=\mathrm{K}_\ell \otimes \mathrm{M}_\ell
					+\mathrm{M}_\ell \otimes \mathrm{K}_\ell$
				and $\mathcal{M}_{\ell}=\mathrm{M}_\ell \otimes \mathrm{M}_\ell$
				are the corresponding stiffness and mass matrices on the face. For $d=3$, we obtain
				\[
					L_{\ell,T}^*=
					 \frac{h_\ell}p \left(
							\mathcal{K}_{\ell}
                      					+  \frac  {p^2} {h_\ell^2} \mathcal{M}_{\ell}  \right) \ .
				\]
				Here, analogously to the case of the patch-interiors, the subspace corrected mass smoother is used.
				Note that the subspace corrected mass smoother is set up such that it bounds the stiffness matrix $\mathcal{K}_{\ell}$
				from above, cf.~\cite[eq.~(11)]{HT:2016}.	In the present paper, besides a trivial scaling, the
				stiffness matrix $\mathcal{K}_{\ell}$ is augmented by $p^2 h_\ell^{-2}$ times
				the mass matrix $\mathcal{M}_{\ell}$. So, we have also to augment the local contributions for the subspace corrected
				mass smoother, cf. the matrices $L_{\alpha}$ in~\cite[Sec.~4.2]{HT:2016}, in the same way.
\end{itemize}

\section{The parallelization of the multigrid solver}\label{sec:parallel}

The parallelization of the multigrid solver follows the approach presented
in~\cite{DHL}. We use MPI\footnote{Message Passing Interface, see \url{http://mpi-forum.org/}.},
so each processor executes independently
the whole algorithm with its local data until communication is explicitly requested.

\begin{figure}[h]
\begin{center}
	\includegraphics[width=.17\textwidth]{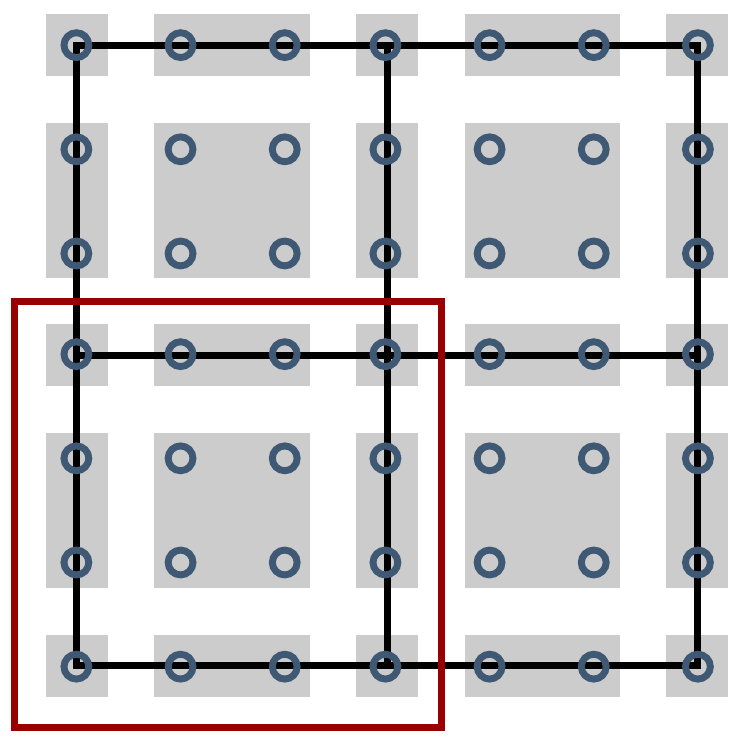}
\end{center}
\caption{The distribution of the dofs to the processors}\label{fig:0b}
\end{figure}

We assign each of the patches to one of the processors. So, that
processor holds the values of all dofs that belong
to that patch including its interfaces, cf. Fig.~\ref{fig:0b}. This means that the dofs
on the interfaces might be assigned to more than one processor.

A vector that occurs in the algorithm, say $\ul{w}$, is stored either in
accumulated form (Type~I) or distributed form (Type~II).
We say that a vector is stored in accumulated form if each of
the processors holds those parts of the global vector which correspond
to the dofs assigned to the processor.
We denote such vectors by $\acc{w}$.
We say that a vector is stored in distributed form if the
global vector is the sum of the contributions of all the processors. Such
vectors are denoted by $\dist{w}$.
Again, the processor-local contributions of the distributed vectors are
supported only the patches assigned to the processor including their
interfaces.

Note that only certain kinds of operations make sense; so we can add
accumulated and distributed vectors only to vectors of the same type:
\begin{equation}\nonumber
\acc{u} + \acc{v}   \rightarrow \acc{w} 
\quad \mbox{and} \quad
\dist{u} + \dist{v}   \rightarrow \dist{w} \ ,
\end{equation}
cf.~\cite[Sec.~5.3]{DHL}.
As the multi-patch setting is equivalent to a standard approach of non-overlapping domain
decomposition, the overall stiffness matrix is assembled on a per-patch basis, i.e., the bilinear
forms $a_k$ from~\eqref{eq:model:bil} are evaluated separately yielding matrices $A_{\ell,k}$. 
Consequently, the global stiffness matrix $A_\ell$ is the sum of
the local contributions $A_{\ell,k}$. This means that the matrix $A_\ell$ is stored in
distributed form, which yields the following mapping type:
\begin{equation}\nonumber
		A_\ell \,  \acc{u} \rightarrow \dist{w}\ ,
\end{equation}
i.e., $A_\ell$ can be applied to accumulated vectors and the the result of the operation is distributed, cf.~\cite[Sec.~5.4.1]{DHL}.

Similar to~\cite[Sec.~7.2.2]{DHL}, the inter-grid transfer operators satisfy 
\begin{equation}\nonumber
P_\ell \, \acc{u}  \rightarrow \acc{w} 
\quad \mbox{and} \quad
P_\ell^\top \,  \dist{u}  \rightarrow \dist{w}
\end{equation}
because the prolongation operator has a block-triangular structure as
in~\cite[eq.~(5.9)]{DHL} and the restriction operator has a block-triangular structure
as in~\cite[eq.~(5.10)]{DHL}. The block-triangular structure is obtained because the
following statements hold true:
\begin{itemize}
	\item On each vertex, the prolonged value $\acc{w}$ coincides with the coarse-grid value $\acc{u}$ of the same vertex.
	\item On each edge, the prolonged values $\acc{w}$ only depend on the coarse-grid values $\acc{u}$ on the same edge and
			on the adjacent vertices.
	\item On each patch-interior, the prolonged values $\acc{w}$ only depend on the coarse-grid values $\acc{u}$ on the same patch-interior and
			on the adjacent edges and vertices.
\end{itemize}
For three dimensions, completely analogous statements hold true.

The global operator $L_\ell^{-1}$ is block-diagonal, where each block corresponds to one piece. Note that
by construction each piece belongs as a whole to one processor or is shared as a whole by
the same processors, so it satisfies both the conditions of~\cite[eq.~(5.9)]{DHL}
and \cite[eq.~(5.10)]{DHL}. This shows
\begin{equation}\nonumber
L_\ell^{-1} \,  \acc{u} \rightarrow \acc{w} 
\quad \mbox{and} \quad
L_\ell^{-1} \,  \dist{u} \rightarrow \dist{w} \ ,
\end{equation}
i.e., this operator can be applied both to distributed and accumulated vectors and it preserves
the type of the vector.

As in any iterative solver, we need to accumulate the vectors of interest in each iterate.
This we denote using the symbol $\Acc$, which maps as follows: 
\begin{equation}\label{eq:ops6}
		\Acc \,  \dist{u}   \rightarrow \acc{w}.
\end{equation}
We note that only a communication between the processors holding neighboring patches is required
in order to perform~\eqref{eq:ops6}.

Only the coarsest grid level $\ell=0$ needs some special treatment. Since the focus of the present
paper is set on parallelizing the multigrid solver without changing its mathematical meaning, we
perform an exact global solve on the coarsest grid level. 
This seems to be acceptable as it is done only for the coarsest
grid level. So, we are required to
communicate the stiffness matrix between all processors such that every
processor holds a global stiffness matrix. 
We set up a direct solver $A_{0}^{-1}$ for this global stiffness matrix, so its application is perform
in the following way
\begin{equation}\nonumber
		\Split A_{0}^{-1} \Accglob \, \acc{u} \rightarrow \acc{w} \ ,
\end{equation}
where $\Accglob$ denotes the accumulation of vectors where each processor obtains the global
vector and $\Split$ is the restriction of the global vector to the patches assigned to
the processor. The latter involves only discarding unnecessary data. We obtain
\begin{align*}
	\Accglob   \dist{u} \to \glob{w} \quad \text{and} \quad \Split \glob{u}\to \acc{w}\ .
\end{align*}

Overall, the parallel multigrid solver looks as follows:

\begin{algorithm}[H]
    \caption{Parallel multigrid solver}
    \label{alg-mg}
    \begin{algorithmic}[1] 
        \Procedure{Multigrid}{$\ell, \acc{u}, \dist{f}$}   
	         \ForAll{$i=1,\ldots,\nu$} \Comment{Pre-smoothing}
		   \State $\acc{u} \gets \acc{u} + \tau \Acc L_\ell^{-1}  ( \dist{f} - A_\ell \acc{u}) $
		\EndFor
		\State $\dist{r} \gets P_\ell^\top ( \dist{f} - A_\ell \acc{u}) $
		\If {$\ell =1$} \Comment{Coarse-grid correction} 
	           	 \State $\acc{p} \gets \Split A_{0}^{-1}\Accglob \dist{r} $ \Comment{Exact solver for coarsest grid level}
		\Else
			\State $\acc{p} \gets 0$ 				
		         	\ForAll{$i=1,\ldots,\mu$}					\Comment{$\mu=1$ is V-cycle; $\mu=2$ is W-cycle}
			   	\State $\acc{p}\gets \textsc{Multigrid}(\ell-1,\acc{p}, \dist{r})$
			\EndFor
	          \EndIf
		\State $\acc{u} \gets \acc{u} + P_\ell \acc{p} $
	         \ForAll{$i=1,\ldots,\nu$} \Comment{Post-smoothing}
		   \State $\acc{u} \gets \acc{u} + \tau \Acc L_\ell^{-1}  ( \dist{f} - A_\ell \acc{u}) $
		\EndFor
		\State \textbf{return } $\acc{u}$
        \EndProcedure
    \end{algorithmic}
\end{algorithm}

We use our multigrid algorithm as a preconditioner for a standard parallel preconditioned conjugate gradient (PCG)
solver. Note that the multigrid preconditioner already takes a distributed residual and returns an accumulated update. So, the
preconditioned conjugate gradient solver only needs to accumulate data in order to compute the required scalar
products accordingly, cf.~\cite[Sec.~6.3.1]{DHL}.

\section{Numerical experiments}\label{sec:num}

In this section, we present numerical experiments concerning the parallelization of the multigrid
solver. The solver was implemented in C++ based on the G+Smo
library~\cite{gismoweb} and, as already mentioned, the parallelization is performed using MPI.
All numerical experiments have been done using the HPC Cluster RADON1\footnote{We use up to 32 out
of 68 available nodes, each equipped with 2x Xeon E5-2630v3 ``Haswell'' CPU (8 cores, 2.4 Ghz, 20 MB cache)
and 128 GB RAM. More information is available at \url{https://www.ricam.oeaw.ac.at/hpc/}.}.

We present timings for setup, assembling and solving.
The setup costs include  
\begin{enumerate}
	\item the costs of the setup of the dof-mappers, which describe the
relation between the local dof-indices and the global dof-indices,
	\item the costs of the grid refinement and the setup of the inter-grid transfer matrices,
	\item the costs of the setup of the piece-local smoothers and
	\item the costs of the setup of the coarse-grid solver.
\end{enumerate}
Here, our implementation of item~1 requires that each processor knows about the indexing
of the global dofs. Also for item~4, the information on all dofs is required, however only
on the coarsest grid level. The costs which are typically dominant, i.e., those for assembling
and for solving, are presented separately. It is important to note that assembling does
not require the any kind of communication between the processors. So its parallelization
is trivial. The communication, which is required for the solving phase,
is discussed in detail in Sec.~\ref{sec:parallel}.

\begin{figure}[h]
\begin{center}
\subfloat[The Yeti footprint\label{fig:1}]{
    \begin{minipage}{.45\textwidth}
        \centering\includegraphics[height=.45\textwidth]{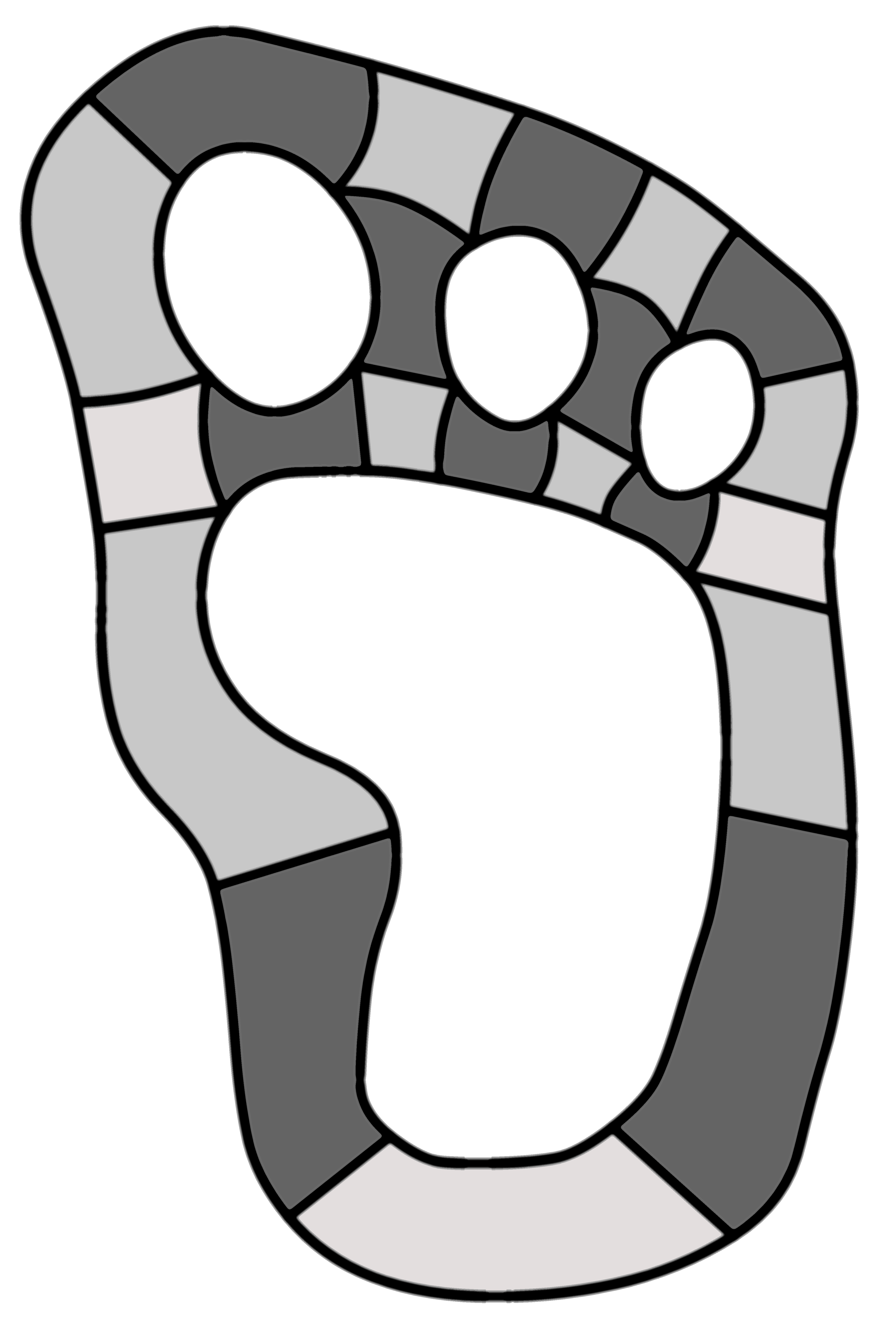}
    \end{minipage}
}
\subfloat[The Fichera corner\label{fig:2}]{
    \begin{minipage}{.45\textwidth}
        \centering\includegraphics[height=.45\textwidth]{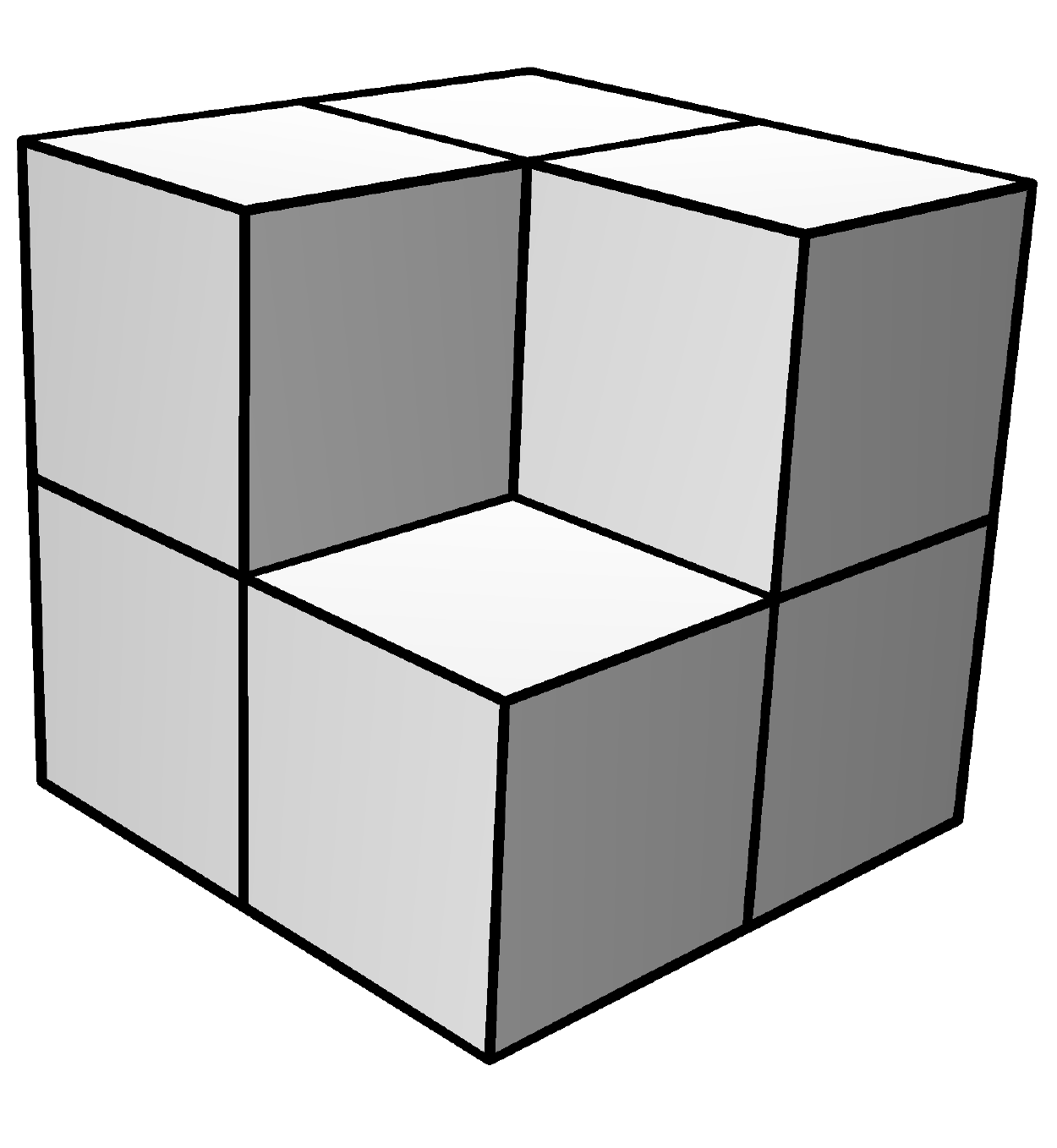}
    \end{minipage}
}
\end{center}
\caption{The computational domains\label{fig:12}}
\end{figure}

We have performed the numerical experiments for two and three dimensions. As
two dimensional domain, we use the Yeti footprint (Fig.~\ref{fig:1}), which has already been
considered in~\cite{Takacs:2017} and which is also a popular domain for the IETI-DP method,
cf.~\cite{KleissEtAl:2012}. As three dimensional domain, we consider the Fichera corner (Fig.~\ref{fig:2}).
This domain
is often considered as extension of the L-shaped domain to three dimensions; the corresponding 
numerical experiments show that the proposed method can also be applied to domains without
full elliptic regularity.

\subsection{The Yeti footprint (2D)}

On the Yeti footprint, we solve the model problem
\[
    \begin{aligned}
        -\Delta u & =50\pi^2\sin(5\pi\ x) \sin(5\pi\ y ) && \quad \mbox{ in } \Omega \ ,\\
                        u & =0&& \quad \mbox{ on } \Gamma_D \ ,\\
            \frac{\partial}{\partial n} u & = 0 && \quad \mbox{ on } \Gamma_N \ ,
    \end{aligned}
\]
where $\Gamma_D$ is the outer boundary and $\Gamma_N$ are the four inner boundaries.

The Yeti footprint consists of $21$ patches, which can be seen in Fig.~\ref{fig:1}. Since we need sufficiently
many patches for parallelization, we first split each patch uniformly into $16$ patches, so we
obtain in total~$K=336$ patches. We solve the problem with a conjugate gradient solver, preconditioned
with one V-cycle of the multigrid method. We perform 1+1 smoothing
steps of the proposed smoother. The damping parameter and the scaling parameter (in the subspace
corrected mass smoother) are chosen as in~\cite{Takacs:2017}, i.e., $\tau = 0.25$ and $\sigma = \tfrac{1}{0.2} h_\ell^{-2}$.

In Tab.~\ref{tab:I:it}, we report on the number of iterations required to reach the desired relative accuracy goal
of $10^{-8}$. Here, $\ell$ represents the number of refinement levels and $p$ the spline degree. On the coarsest
grid level ($\ell=0$), the patch-local discretization only consists of global polynomials, i.e., each patch is one
element of the discretization. Refinement is done by uniformly refining the patch-local grids, keeping the number
of patches unchanged. We observe, as in~\cite{Takacs:2017}, that the number of iterates is quite robust in the
grid size and in the spline degree. The presented numbers have been computed with the serial code. The number of
iterates is supposed to be the same if parallelization is applied; however due to some small numerical
instabilities, in some cases the parallel code needs one additional iteration (but never more than that).
Similar iteration counts are obtained for the W-cycle.

In Tab.~\ref{tab:I:strong}, we present the strong scaling results. We fix the grid level $\ell$ and the spline
degree $p$ to two typical values. For $\ell=7$ and $p=4$, we have $5\,768\,189$ dofs and the corresponding
stiffness matrix has $4.6\;10^8$ non-zero entries. For the case $\ell=7$ and $p=8$, the number of dofs
increases slightly to $6\,125\,757$, but the stiffness matrix has already $1.8\;10^9$ non-zero entries.
In the first two rows, we compare the costs of the
serial code and the parallel code. Here, we obtain that the parallel code is slightly slower during the
solving phase which is mainly due to the fact that the parallel code does not assemble the whole stiffness
matrix but works with patch-local stiffness matrices. This allows also to consider the larger problem
with $\ell=7$ and $p=8$, where the serial code caused memory problems.

\begin{table}[t]
\begin{center}
    \begin{tabular}{l|rrrrrrr}
    \toprule
    \multicolumn{1}{l}{$\ell\;\diagdown\; p$}    & \quad 2 & \quad 3 & \quad 4 & \quad 5 &\quad  6 &\quad  7 &\quad  8  \\
    \midrule
    4                 & 45 & 42 & 37 & 33 & 31 & 28 & 25  \\
    5                 & 48 & 44 & 40 & 36 & 33 & 30 & 27  \\
    6                 & 50 & 44 & 41 & 36 & 35 & 33 & 27  \\
    7                 & 51 & 45 & 42 & 37 & 36 & 34 & 28  \\                                                       
    \bottomrule                                                                                               
    \end{tabular}                                                                                             
\end{center}                                                                                                      
\caption{Iteration counts for Yeti footprint, $K=336$}                                   
\label{tab:I:it}                                                                                                  
\end{table}                                                                                                       

\begin{table}[t]
\begin{center}
    \begin{tabular}{r||rr|rr|rr||rr|rr|rr}
    \toprule
    & \multicolumn{6}{c||}{$\ell= 7\ ,\;\;p=4\ ,\;\;K=336$} & \multicolumn{6}{c}{$\ell= 7\ ,\;\;p=8\ ,\;\;K=336$} \\
    \midrule
                 & \multicolumn{2}{l|}{Setup} &  \multicolumn{2}{l|}{Assembling}  & \multicolumn{2}{l||}{Solving} &
                   \multicolumn{2}{l|}{Setup} &  \multicolumn{2}{l|}{Assembling}  & \multicolumn{2}{l}{Solving} \\
     $\#$  Proc. &$t$&$s$&$t$&$s$&$t$&$s$&$t$&$s$&$t$&$s$&$t$&$s$  \\
    \midrule 
        (serial)   & 217.1  & \; --\,\,  & 522.1       & \; --\,\, & 4929.5         & \;  --\,\,   & --\,\,     &   \; --\,\,& --\,\,     & \;  --\,\,    & --\,\,      & \;  --\,\,   \\
         1         & 220.0  & \; 1\;\;\, & 520.0       & \; 1\;\;\,& 5125.9         & \;  1\;\;\,  & 549.0      &   \;1\;\;\,& 8230.9     & \;   1\;\;\,  & 4729.8      & \;  1\;\;\,  \\
         2         & 80.1   & \;   2.7   & 263.8       & \;   1.9  & 1367.6         & \;   3.7     & 225.4      &   \;  2.4  & 4158.5     & \;     1.9    & 1250.8      & \;    3.7    \\
         4         & 35.3   & \;   6.2   & 131.8       & \;   3.9  & 399.1          & \;  12.8     & 117.0      &   \;  4.6  & 2098.5     & \;     3.9    & 409.9       & \;   11.5    \\
         8         & 17.1   & \;  12.8   & 66.0        & \;   7.8  & 109.6          & \;  46.7     & 53.9       &   \; 10.1  & 1055.9     & \;     7.8    & 140.2       & \;   33.7    \\
         16        & 10.7   & \;  20.5   & 33.9        & \;  15.3  & 40.7           & \; 125.9     & 30.0       &   \; 18.3  & 543.4      & \;    15.1    & 59.2        & \;   79.9    \\
         32        & 8.0    & \;  27.5   & 17.4        & \;  29.8  & 17.1           & \; 299.7     & 17.7       &   \; 31.0  & 275.1      & \;    29.9    & 26.8        & \;  176.4    \\
         64        & 7.2    & \;  30.5   & 9.4         & \;  55.3  & 10.6           & \; 483.5     & 12.9       &   \; 42.5  & 149.7      & \;    54.9    & 13.7        & \;  345.2    \\
         128       & 6.0    & \;  36.3   & 5.1         & \; 101.3  & 4.1            & \;1250.2     & 9.9        &   \; 55.4  & 76.2       & \;   108.0    & 7.2         & \;  656.9    \\
         256       & 6.3    & \;  34.4   & 3.2         & \; 160.0  & 3.3            & \;1553.3     & 9.5        &   \; 57.7  & 51.4       & \;   160.1    & 6.5         & \;  727.6    \\
    \bottomrule                                                                                                      
    \end{tabular}
\end{center}
\caption{Strong scaling behavior for Yeti footprint}
\label{tab:I:strong}
\end{table}

\begin{table}[t]
\begin{center}
    \begin{tabular}{r||r|r|r|r|r||r|r|r|r|r}
    \toprule
    & \multicolumn{5}{c||}{$\ell= 7\ ,\;\;p=4$} & \multicolumn{5}{c}{$\ell= 7\ ,\;\;p=8$} \\
	\midrule
    $\#$  Proc. & \multicolumn{1}{p{2.5em}|}{$\#$ dofs}  & \multicolumn{1}{p{1.5em}|}{It.}&    \multicolumn{1}{p{3em}|}{Setup} &  \multicolumn{1}{p{3.5em}|}{Ass.}  &  \multicolumn{1}{p{3.5em}||}{Solving}  
                & \multicolumn{1}{p{2.5em}|}{$\#$ dofs}  & \multicolumn{1}{p{1.5em}|}{It.}&    \multicolumn{1}{p{3em}|}{Setup} &  \multicolumn{1}{p{3.5em}|}{Ass.}  &  \multicolumn{1}{p{3.5em}}{Solving}   \\
    \midrule
         4         & 360\,902  & 46     & 1.4     &  9.2       & 7.9        & 383\,262     & 46   & 6.6   &   147.3         & 21.4    \\
         16        & 1\,442\,569  & 44  & 2.4     &  9.3       & 8.5        & 1\,531\,977  & 44   & 8.7    &  151.2         & 20.8    \\
         64        & 5\,768\,189  & 41  & 7.2     &  9.5       & 10.5        & 6\,125\,757  & 28  & 13.2    & 148.7        & 13.7        \\   
         256       &23\,068\,573  & 36  & 42.4    &  9.7       & 9.5         &24\,498\,717  & 26  & 54.5    &  153.8        & 20.0       \\     
    \bottomrule                                                                                                      
    \end{tabular}                                                                 
\end{center}                                                                                          
\caption{Weak scaling behavior for Yeti footprint}                                  
\label{tab:I:weak}                                                                  
\end{table}

In the following rows, we consider the strong scaling behavior. We present in each case the time $t$ in
seconds required for setup, assembling and solving and the corresponding speedup $s$.
We observe that the overall method has good strong scaling properties. As the setup phase consists also
of parts that are not parallelized, we observe this time does not fall below a few seconds. The assembling
phase, which is known to be dominant phase in high-order isogeometric methods, scales almost optimal. Also the
solving phase needs rather little communication and is expected to scale well therefore. Indeed, the speedup
is much larger than what would be expected. The authors think that this might be explained by some extraordinary
caching effects, but here further investigation is required. The extraordinary well behavior of the solver
cannot be explained with changed convergence behavior because in all cases, the convergence behavior is
identical.

In Tab.~\ref{tab:I:weak}, we present weak scaling results. We again fix the grid level $\ell$ and the spline
degree $p$ to two typical values. Here, for the case of $4$ processors, we consider the initial configuration of
$K=21$ patches; in the following rows we consider $84$, $336$ and $1344$ patches. (So, the third row
with $64$ processors coincides with the line with $64$ processors in Tab.~\ref{tab:I:strong}.)
As the setup phase is not fully parallelized, the setup times increase if the number of patches is increased.
Both, the assembling times and the solving times are rather constant and do not indicate a clear tendency.
The solving times also change due to the fact that the required number of iterations decays if the patches
are split up. The computational costs for the global-coarse grid solver is negligible in this example;
for $K=1344$ patches the costs are $0.36$ seconds for $p=4$ and $1.4$ seconds for $p=7$ and for smaller
patch numbers even less.

\subsection{The Fichera corner (3D)}

On the Fichera corner, we solve the model problem
\[
    \begin{aligned}
        -\Delta u & =  75\pi^2\sin(5\pi\ x) \sin(5\pi\ y)\sin(5\pi\ z) &&  \mbox{ in } \Omega := (0,2)^3 \backslash [1,2)^3 \ ,\\
                        u & =0 &&  \mbox{ on } \Gamma_D := \{ (x,y,z) \in \partial \Omega \;:\; xyz=0 \} \ ,\\
            \frac{\partial}{\partial n} u & = 0 &&  \mbox{ on } \Gamma_N := \partial \Omega \backslash \Gamma_D \ .
    \end{aligned}
\]

The Fichera corner consists of $7$ patches, which can be seen in Fig.~\ref{fig:2}, which are uniformly
split into~$K=448$ patches in total. Again, we solve the problem with a conjugate gradient solver preconditioned
with one V-cycle of the multigrid method with 1+1 smoothing steps. Again
$\tau = 0.25$ and $\sigma = \tfrac{1}{0.2} h_\ell^{-2}$ are chosen.

In Tab.~\ref{tab:II:it}, we report on the number of iterations required to reach the desired relative
accuracy goal of $10^{-8}$.  We observe, as for the Yeti footprint, that the number of iterates is
quite robust in the grid size and in the spline degree. The presented numbers have been computed with
the serial code. Again, the parallel code yields (almost) the same numbers.

In Tab.~\ref{tab:II:strong}, we present the strong scaling results. 
For $\ell=4$ and $p=2$, we have $N=2\,201\,024$ dofs and a stiffness matrix with $2.5\;10^8$ non-zero entries.
The second example with $\ell=3$ and $p=4$ yields $N=596\,288$ dofs and $2.8\;10^8$ non-zero entries.
The timings behave similar as in the two-dimensional case, however the costs of the setup phase are
much larger which can be explained by the fact that the interfaces are much larger.
(For two dimensional problems, the interfaces consist of
$\mathcal{O}(N^{1/2})$ dofs and for three dimensional problems, the interfaces consist of
$\mathcal{O}(N^{2/3})$ dofs.) The assembling times seem to be optimal, whereas the solving times
again behave extraordinary well.

\begin{table}[t]
\begin{center}
    \begin{tabular}{l|rrrrrr}
    \toprule
    \multicolumn{1}{l}{$\ell\;\diagdown\; p$}    & \quad 2 & \quad 3 & \quad 4 & \quad 5 &\quad  6  \\
    \midrule
    1                 & 30 & 31 & 31 & 26 & 22   \\
    2                 & 33 & 32 & 33 & 31 & 28   \\
    3                 & 39 & 38 & 37 & 33 & 30   \\
    4                 & 44 & 44 & 42 & 37 & 35   \\
	\bottomrule
	\end{tabular}                                                                                                          
\end{center}                                                                                                                   
\caption{Iteration counts for Fichera corner, $K=448$}                                                                                  
\label{tab:II:it}                                                                                                              
\end{table}

\begin{table}[t]                                                                                                               
\begin{center}
    \begin{tabular}{r||rr|rr|rr||rr|rr|rr}
    \toprule
    & \multicolumn{6}{c||}{$\ell= 4\ ,\;\;p=2\ ,\;\;K=448$} & \multicolumn{6}{c}{$\ell= 3\ ,\;\;p=4\ ,\;\;K=448$} \\
    \midrule
                 & \multicolumn{2}{l|}{Setup} &  \multicolumn{2}{l|}{Assembling}  & \multicolumn{2}{l||}{Solving} &
                   \multicolumn{2}{l|}{Setup} &  \multicolumn{2}{l|}{Assembling}  & \multicolumn{2}{l}{Solving} \\
     $\#$  Proc. &$t$&$s$&$t$&$s$&$t$&$s$&$t$&$s$&$t$&$s$&$t$&$s$  \\
    \midrule
    (serial)               &   179.4   & \; --\,\, & 260.2       & \; --\,\, & 4980.7  & \;  --\,\,   &      93.5    & \; --\,\,&  1313.5     & \; --\,\,  &   1252.2      & \;--\,\,   \\
         1                 &   198.2   & \; 1\;\;\,& 253.4       & \; 1\;\;\,& 5091.3  & \;    1\;\;\,&     109.7    & \;1\;\;\,&  1985.9     & \; 1\;\;\, &   1073.1      & \; 1\;\;\, \\
         2                 &   77.7    & \; 2.5    & 127.4       & \;   1.9  & 1355.0  & \;     3.7   &      51.2    & \;  2.1  &  1103.3     & \;   1.8   &    340.0      & \;    3.1  \\
         4                 &   49.2    & \; 4.0    & 63.5        & \;   3.9  & 395.7   & \;    12.8   &      30.6    & \;  3.5  &  492.2      & \;   4.0   &    110.2      & \;    9.7  \\
         8                 &   32.9    & \; 6.0    & 32.0        & \;   7.9  & 99.2    & \;    51.3   &      22.8    & \;  4.8  &  214.1      & \;   9.2   &     40.1      & \;   26.7  \\
         16                &   28.3    & \; 7.0    & 16.5        & \;  15.3  & 34.4    & \;   148.0   &      27.3    & \;  4.0  &  113.4      & \;  17.5   &     22.7      & \;   47.2  \\
         32                &   26.8    & \; 7.4    & 8.4         & \;  30.1  & 12.4    & \;   410.5   &      17.0    & \;  6.4  &  55.5       & \;  35.7   &      8.3      & \;  129.2  \\
         64                &   32.2    & \; 6.1    & 5.5         & \;  46.0  & 5.0     & \;  1018.2   &      24.4    & \;  4.5  &  31.2       & \;  63.6   &      6.9      & \;  155.5  \\
         128               &   42.3    & \; 4.6    & 2.7         & \;  93.8  & 2.6     & \;  1958.1   &      15.8    & \;  6.9  &  15.4       & \; 128.9   &      3.6      & \;  298.0  \\
         256               &   54.1    & \; 3.6    & 1.1         & \; 230.3  & 1.8     & \;  2828.5   &      24.0    & \;  4.5  &  7.9        & \; 251.3   &      3.9      & \;  275.1  \\
	\bottomrule                                                                                                      
	\end{tabular}
\end{center}
\caption{Strong scaling behavior for Fichera corner}
\label{tab:II:strong}
\end{table}

\begin{table}[t]
\begin{center}
    \begin{tabular}{r||r|r|r|r|r||r|r|r|r|r}
    \toprule
    & \multicolumn{5}{c||}{$\ell= 4\ ,\;\;p=2$} & \multicolumn{5}{c}{$\ell= 3\ ,\;\;p=4$} \\
    \midrule
    $\#$  Proc. & \multicolumn{1}{p{2.5em}|}{$\#$ dofs}  & \multicolumn{1}{p{1.5em}|}{It.}&    \multicolumn{1}{p{3em}|}{Setup} &  \multicolumn{1}{p{3.5em}|}{Ass.}  &  \multicolumn{1}{p{3.5em}||}{Solving}  
                & \multicolumn{1}{p{2.5em}|}{$\#$ dofs}  & \multicolumn{1}{p{1.5em}|}{It.}&    \multicolumn{1}{p{3em}|}{Setup} &  \multicolumn{1}{p{3.5em}|}{Ass.}  &  \multicolumn{1}{p{3.5em}}{Solving}   \\
    \midrule                                                                           
         1         &  34\,391     & 28 & 0.4   & 4.0         & 1.9     &  9\,317      & 31 &  0.6   &    38.0         &    1.3            \\
         8         & 275\,128     & 39 & 1.3   & 4.5         & 4.0     & 74\,536      & 35 &  1.1   &    22.8         &    2.1            \\
         64        &2\,201\,024   & 45 & 54.4  & 4.5         & 6.6     &596\,288      & 38 &  24.1  &    28.7         &    6.2            \\     
         512       &17\,608\,192  & 46 &2071.3 & 5.1         &  11.1   &4\,770\,304   & 35 &  2343.8 &   32.4         &   59.5           \\       
    \bottomrule                                                                                                      
    \end{tabular}                                                                                                   
\end{center}
\caption{Weak scaling behavior for Fichera corner}
\label{tab:II:weak}
\end{table}

In Tab.~\ref{tab:II:weak}, we present weak scaling results. We again fix the grid level $\ell$ and the spline
degree $p$ to two typical values. Here, for the case of $4$ processors, we consider the initial configuration
of $K=7$ patches. For the following rows, we consider $56$, $448$ and $3584$ patches.
(So, the line with $64$ processors coincides with the corresponding line in Tab.~\ref{tab:II:strong}.)
Again, the assembling times and the solving times do not show
any clear tendency. Only for the last line with $3584$ patches, the coarse-grid
solver causes problems. For the case $\ell=3$ and $p=4$, $51$ of the $59$ seconds required for solving are
due to the global solver on the coarsest grid. Again, the setup costs get dominant if the number of patches
is increased.

Concluding, we have shown that the robust multi-patch multigrid solver
from~\cite{Takacs:2017} can be extended to three dimensional domains and that it converges
well also in this case. We have observed that the multigrid solver can be
parallelized in a natural way yielding very good speedup rates.
Certainly, this is not the end of the story and further improvement should
be considered in two directions. First, the setup phase becomes a bottleneck if many
processors are considered. Here, improvements would be mainly a challenge in terms of
implementation and data management. Second, the coarse-grid problem becomes too large if the
number of patches is increased, particularly in the three dimensional case. To resolve that
issue, it would be necessary to further coarsen the coarse-grid problem or to consider
approximate solvers on the coarsest grid level which certainly would change the mathematical
meaning of the algorithm and could, therefore, influence its
convergence behavior. Finally, further investigation is required
to completely understand the super optimal speedup rates observed in the
strong scaling tests.\\[.25em]

\textbf{Acknowledgments.}
The first author would like to thank the Austrian Science Fund (FWF)
for the financial support through the DK W1214-04, while the second
author was supported by the FWF grant NFN S117-03.

\bibliographystyle{amsplain}
\bibliography{references}

\end{document}